\documentstyle[12pt]{article}

\font\tenmsb=msbm10
\font\sevenmsb=msbm7
\font\fivemsb=msbm5

\newfam\msbfam
\textfont\msbfam=\tenmsb
\scriptfont\msbfam=\sevenmsb
\scriptscriptfont\msbfam=\fivemsb
\def\Bbb#1{{\fam\msbfam #1}}

\makeatletter
\ifnum\@ptsize=0  \fi \ifnum\@ptsize=2 
 \fi 

\newcommand\qed{{\hspace*{\fill}Q.E.D.\vskip12pt plus 1pt}} 

\newcommand\sK{{\cal K}}

\newcommand\sO{{\cal O}}

\newcommand\bn{{\bigskip \noindent}}

\newcommand\sn{{\smallskip \noindent}}

\newcommand\zed{{\Bbb Z}}

\newcommand\bQ{{\Bbb Q}}
\newcommand\bR{{\Bbb R}}

\newcommand\bN{{\Bbb N}}
\newcommand\Int{{\rm Int}}

\newcommand\cK{{\overline {{\cal K}(X)} \setminus \{0\}}}

\newcommand\bP{{\Bbb P}}

\newcommand\proof{{\noindent\bf Proof.\ }} 

\newtheorem{theorem}{Theorem}[section]

\newtheorem{corollary}[theorem]{Corollary} 
\newtheorem{proposition}[theorem]{Proposition} 
 
\newtheorem{re}[theorem]{Remark}

\newtheorem{definition}[theorem]{Definition} 
\newtheorem{conjecture}[theorem]{Conjecture} 
 
\newtheorem{problem}[theorem]{Problem}
\newtheorem{example}[theorem]{Example}
\newtheorem{notations}[theorem]{Notations}

\newenvironment{remark}{\begin{re}\em}{\end{re}}

\begin{document}
\title { The dual K\"ahler cone of compact K\"ahler threefolds } \author{ Keiji Oguiso \and 
Thomas Peternell}

\maketitle

\vspace*{-0.5in}\section*{Introduction} The K\"ahler cone and its dual play an important role in the 
study of compact
K\"ahler manifolds. Therefore it seems natural to ask whether one can read off the (dual) 
K\"ahler cone, whether the
underlying manifold is projective or not. A classical theorem of Kodaira says that a compact 
K\"ahler manifold whose K\"ahler
cone has an interior rational point, is projective. Indeed, a multiple of such a point defines 
a Hodge metric on the 
manifold. But not only the K\"ahler cone itself, also its dual in $H^{n-1,n-1}(X)$ (1.7), 
where $n$ is the dimension of the K\"ahler
manifold $X$, is of interest. This is parallel to the projective theory, where 
both the ample cone and its dual cone ${\overline {NE}}$,
the Mori cone of curves,
play an important role. One of the basic questions underlying this paper asks whether 
projectivity of K\"ahler manifolds
can be expressed in terms of the dual K\"ahler cone. This question was first posed and treated for surfaces by Huybrechts [Hu99];
an geometric proof for surfaces was given in [OP00]. To be precise, we consider the following
\begin{problem} Let $X$ be a compact K\"ahler manifold of dimension $n$ and suppose that the dual K\"ahler cone $\sK^*(X) \subset H^{n-1,n-1}(X)$ contains a rational
interior point. How algebraic is $X$, i.e. what can be said about the algebraic dimension $a(X)$ of $X$?
\end{problem}
Since this problems seems to be rather hard, we restrict ourselves to special rational points, namely those
represented by {\it effective} curves. Here is a rather partial answer to the problem.

\begin{theorem} Let $X$ be a smooth compact K\"ahler threefold, $C \subset X$ an effective curve such that $[C]Ê\in \Int \sK^*(X),$ the interior of 
the dual K\"ahler cone.
Then $a(X) \geq 2$ unless $X$ is simple non-Kummer.
\end{theorem}

The proof relies on classification theory of compact K\"ahler threefolds; the term ``simple" means that $X$ does not admit a covering
family of curves. It is expected that such a simple threefold should be Kummer, i.e. bimeromorphic to a quotient of a torus by a
finite group. This conjecture is however still far from being proven, therefore we need to make the exception in the theorem.
It seems not impossible that Theorem (0.2) is actually sharp; in (3.2)-(3.4) we give a potential procedure how to construct a K\"ahler threefold with $a(X) = 2$
and admitting an irreducible curve which represents an interior rational point of the dual K\"ahler cone. \\
Given a curve $C \subset X$, we would like to decide whether its class $[C]$ is an interior point of $\sK^*(X)$ or a boundary point. Obvious examples for boundary
points are provided by curves contracted by maps $\phi: X \to Y$ to K\"ahler spaces. It seems therefore very reasonable to expect the following

\begin{conjecture} If the normal bundle $N_{C/X}$ is ample, then $[C] \in \Int \sK^*(X).$
\end{conjecture}

If some multiple of $C$ moves in a family that covers $X$, then $[C] \in \Int \sK^*(X)$
by (4.13). Moreover
we prove (0.3) for complete intersections of hyperplane sections. Combining this with problem (0.1) we are lead to

\begin{problem}Let $C \subset X$ be a smooth curve in the compact K\"ahler manifold $X$ with ample normal bundle. How projective is $X?$
\end{problem}

Our expectation - and we prove this in many cases - is that at least the algebraic dimension $a(X) \geq 2;$ and similarly as in (0.2) there seem to be examples where
$a(X) = 2.$
Of course a similar question can be asked for submanifolds of
larger dimension. We prove

\begin{theorem} Let $C \subset X$ be a smooth curve of genus $g$ with ample normal bundle $N_C$ in the compact K\"ahler manifold $X.$ Then
\begin{description}
\item {(1)} If $g \leq 1$, then $X$ is projective. The same holds for general $g$, if $N_C \otimes T_C$ is ample, where $T_C$ denotes the tangent bundle
of $C.$ 
\item {(2)} If $\kappa (X) \leq 0,$ then $X$ is projective or simple non-Kummer.
\item {(3)} The algebraic dimension $a(X) \ne 1$. 
\item {(4)} If $\dim X = 3$ and $C$ moves in a covering family, then $X$ is projective. 
\end{description}
\end{theorem}

Finally we investigate threefolds containing an elliptic curve $C$ with ample normal bundle. In particular we prove that either $X$ is
rationally connected, so $\pi_1(X) = 0$ or $\pi_1(X) = \zed^2,$ so that 
$$ {\rm Im}(\pi_1(C) \to \pi_1(X)) $$
has always finite index. One might asks whether this holds in general for submanifolds with ample normal bundle.


\tableofcontents

\section{ Preliminaries}

Here we collect basic results and concepts which are essential for this paper. \\
If $X$ is an irreducible reduced compact complex space, we denote by $a(X) $ its algebraic dimension. For this notion, for the notion of an
algebraic reduction and related stuff we refer e.g. to [Ue75],[GPR94]. We say that a meromorphic 
map $f: X \rightharpoonup Y$ is {\it almost holomorphic},
if it is proper and holomorphic on some open non-empty set $U \subset X$. This means that $f$ has 
some compact fibers, i.e. theses fibers do
not meet the set of indeterminacies.

\begin{definition} {\rm Let $X$ be a compact K\"ahler manifold. $X$ is called {\it algebraically 
connected} if there exists a family of curves $(C_t)$
such that $C_t$ is irreducible for general $t$ and such that every two very general points 
can be joined by a chain of $C_t$'s.}
\end{definition}

Then we have Campana's theorem [Ca81]

\begin{theorem}Every algebraically connected compact K\"ahler manifold is projective. 
\end{theorem} 

An immediate consequence of this theorem is

\begin{corollary} Every algebraic reduction of a threefold with $a(X) = 2$ is almost holomorphic. 
\end{corollary}

\begin{definition} {\rm  A compact K\"ahler manifold $X$ is {\it simple}, if there is no proper positive-dimensional subvariety
through a very general point of $X.$ }
\end{definition}

Concerning the structure of simple compact K\"ahler threefolds one has the

\begin{conjecture}  Every simple compact K\"ahler threefold is Kummer, i.e. bimeromorphic to $T/G$, with a torus $T$
and a finite group $G$ acting on $T.$
\end{conjecture}
For the relation to Mori theory in the K\"ahler case, see [Pe98]. 

We will need the following classification of compact K\"ahler threefolds due to Fujiki [Fu83].

\begin{theorem}  Let $X$ be a compact K\"ahler threefold with $a(X) \leq 1. $\\
(1) If $a(X) = 0$ but not simple, then $X$ is uniruled and moreover there exists a holomorphic map $f: X \to S$ to a 
normal surface $S$ with $a(S) = 0,$ such that the general fiber is $\bP_1.$ \\
(2) If $a(X) = 1$ and if a holomorphic model of the algebraic reduction admits a multi-section, then $X$ is bimeromorphic to $(A \times F)/G,$
where $A$ is a torus or K3 with $a(A) = 0,$ where $F$ is a smooth curve and where $G$ is finite group acting on $A$ and on $F$
and acting on the product diagonally. The map $X \rightharpoonup {F/G}$  is an algebraic reduction over the smooth curve ${F/G}.$ \\
If no holomorphic model of the algebraic admits a multi-section, then $f$ is holomorphic and the general fiber is a torus or K3.  

\end{theorem}

\begin{notations} Let $X$ be a compact K\"ahler manifold of dimension $n.$ \\
(1)  The K\"ahler cone $\sK(X)$ is the subset of $H^{1,1} := H^{1,1}(X) \cap H^2(X,\bR)$ 
consisting of the K\"ahler classes of $X.$ It is an open cone in $H^{1,1}.$   \\
(2) The dual K\"ahler cone $\sK^*(X)$ is the dual cone in $H^{n-1,n-1} := H^{n-1,n-1}(X) 
\cap H^{2n-2}(X,\bR)$ 
with respect to the natural pairing $H^{1,1}(X) \times H^{1,1}(X) \to \bR.$  \\
\end{notations}

By ${\overline {\sK(X)}}$ we will denote the closure of the K\"ahler cone in $H^{1,1}$ with respect to the usual topology.
In contrast,the dual K\"ahler cone is closed by definition. \\
The following seems to be well-known, however we could not find an explicit reference, so we give a short proof. A general good 
reference for the theory of current is [Ha77].

\begin{proposition} Let $X$ be a compact K\"ahler manifold. Then $ {\sK^*(X)} $ is the cone ${\cal P}(X)$ of classes of positive closed
currents of bidimension $(1,1).$ 
\end{proposition}

\proof By Demailly-Paun [DP04], Corollary (0.3),  $ {\sK^*(X)} $ is the closed cone generated by classes of currents $[Y] \wedge \omega^{p-1},$
where  $Y$ is an irreducible analytic set of dimension $p$ and $\omega$ a K\"ahler form. Of course, $[Y]$ is the current given by integration
over $Y.$ Since all these currents $[Y] \wedge \omega^{p-1},$ are positive, $ {\sK^*(X)} $ is contained in the cone ${\cal P}(X)$.
The other inclusion being clear, the assertion follows.
\qed

\section{Blow-ups and Galois covers}

We investigate the behaviour of interior points of the dual K\"ahler cone under blow-ups and special Galois covers. The
results will be needed in sect. 3.
In this section $X$ always denotes a compact K\"ahler {\bf threefold} unless otherwise stated.

\begin{proposition} Let $\pi: \hat X \to X$ be the blow-up along a submanifold $Y.$ If $\sK^*(X)$ contains a rational interior point, then
so does $\sK^*(\hat X).$
\end{proposition}

\bn Before giving the proof, we explain the idea in the simple case that the rational interior point,
which is in general represented by a positive closed current, is represented by an effective curve.
So let $C = \sum a_i C_i$ be an effective curve, $a_i > 0$ such that $[C] \in {\rm Int}\sK^*(X).$
In case $C_i \neq Y,$ we let $\hat C_i \subset \hat X$ be the strict transform of $C_i;$ if $C_i = Y,$
we let $\hat C_i$ be a section of $\pi \vert E \to Y.$ \\
Then $\sum a_i [\hat C_i] + m[l]$ is an interior point of  $\sK^*(\hat X)$ for a suitable rational $m.$
Here $l$ is a fiber of $\pi \vert E$ if $\dim Y = 1$ resp. a line in $E \simeq \bP_2$ if $\dim Y = 0.$ \\
In the general case however, the arguments get more involved since ``strict transforms of currents'' cannot be defined in general
[Me96,p.52/53]. This surprising point was explained to us by J.P.Demailly. 
\bn

\proof Let $\alpha \in \Int \sK^*(X) \cap H^4(X,\bQ)$ and represent $\alpha$ by a positive current $T.$ Let $E = \pi^{-1}(Y)$ be the
exceptional divisor and $l$ either a line in $E$, if $\dim Y = 0$ or a ruling line, if $\dim Y = 1.$ We fix a section $\hat Y \subset E$
in case $E$ is ruled. If $E = \bP_2,$ then we set formally $\hat Y = 0.$ \\ 
We consider the canonical decomposition
$$ T = \chi_YT + \chi_{X \setminus Y}T;$$
$\chi_Y$ denoting the characteristic function of $Y.$ All currents occuring in this decomposition are closed. 
If $\dim Y = 0,$ then $\chi_YT = 0;$ if $\dim Y = 1,$ then by Siu's theorem [Si74], $\chi_YT = a T_Y,$ where $a \geq 0$ and $T_Y$ is the
current ``integration over $Y$". In particular  $d\chi_{X \setminus Y}T = d\chi_YT = 0.$ Let $T' := \chi_{X \setminus Y}T $ for
simplicity. We will proceed in three steps.
\begin{enumerate}
\item $\pi^*([T']) \in \sK^*(\hat X)$ (possibly on the boundary);
\item $\hat \alpha := a  [\hat Y] + b[l] + \pi^*([T']) \in {\rm Int} \sK^*(\hat X)$ for $b >  0;$
\item  $\hat \alpha $ can be chosen rational.  
\end{enumerate} 
Starting with the proof of (1), let $\hat \omega$ be a K\"ahler form on $\hat X$, and we need to show
$$ \pi^*([T']) \cdot \hat \omega \geq 0. \eqno (*) $$
Let $S := \pi_*(\hat \omega);$ then $S$ is a positive closed current on $\hat X$ which is smooth outside $Y.$ 
We have
$$ [\hat \omega] = \pi^*[S] + \lambda [E],$$
with some negative number $\lambda.$ 
Therefore 
$$ \pi^*[T'] \cdot \hat \omega = \pi^*[T'] \cdot \pi^*[S] = [T'] \cdot [S]. $$ 
We show that $$ [T'] \cdot [S] \geq 0.$$
In fact, using Demailly's regularization theorem, we can write [S] as a weak limit of smooth positive
closed forms $\Theta_{\epsilon}$ in the same cohomology class as $T'$ such 
$$ \Theta_{\epsilon} \geq - \lambda_{\epsilon}u - O(\epsilon)\eta,$$
where $\eta$ is a positive $(1,1)-$form, $u$ a suitable semi-positive $(1,1)-$form and
$(\lambda_{\epsilon})$ a decreasing family of non-negative smooth functions for 
$0 < \epsilon < 1$ converging pointwise to $0$ on $X \setminus Y$ and to the Lelong number
$\nu(S,x)$ for $x \in Y.$ We conclude
$$ [S] \cdot [T'] = [\Theta_{\epsilon}] \cdot [T'] = T'(\Theta_{\epsilon}) \geq - \int_X \lambda_{\epsilon} u \wedge T' - O(\epsilon).$$
Then the monotone convergence theorem gives the claim since $\chi_YT' = 0.$ 
\\ (2) By (1) we know already that $\hat \alpha \in \sK^*(\hat X)$. If it is not in the interior, then there exists 
$\hat \beta \in \partial \sK(\hat X)$ such that
$$ \hat \alpha \cdot \hat \beta = 0.$$ 
Hence
$$ 0 = a [\hat Y] \cdot \hat \beta + b [l] \cdot \hat \beta + \pi^*([T']) \cdot \hat \beta. \eqno (a)$$
Since all summands are semi-positive, we conclude first
$$ [\hat Y] \cdot \hat \beta = [l] \cdot \hat \beta = 0. \eqno (b) $$ 
Now (b) implies that $\hat \beta = \pi^*(\beta)$, with $\beta \in \sK(X),$
as we shall see in a moment (claim (c)). 
Since by (a) and (b)  we have $[T'] \cdot \beta = 0,$ and since $[T]$ is an interior point of $\sK^*(X),$
we conclude
$$ [Y] \cdot \beta > 0 $$
(this already settles the case $\dim Y = 0$). But then $[\hat Y] \cdot \hat \beta > 0 $ ( represent $\beta$ by a form),
contradiction. \\
It remains to prove claim (c) in order to settle (2).  It is clear by (b)  that $\hat \beta = \pi^*(\beta),$ 
it only remains to show that $\beta \in {\overline { \sK(X)}}.$ Notice that $\beta \cdot Y = 0.$ 
Assuming the contrary, there exists $[T] \in \sK^*(X),$ such that $\beta \cdot [T] < 0.$ By virtue of the decomposition
$$ T = aT_Y + \chi_{X \setminus Y}T $$
we may assume $\chi_Y T = 0.$ Hence by (1)
$$ \hat \beta \cdot \pi^*[T] = \beta \cdot [T] < 0,$$
with $\pi^*[T] \in \sK^*(\hat X),$ contradicting $\hat \beta \in {\overline {\sK(\hat X)}}.$ So (c) and therefore (2) are proved. \\
(3) It remains to show that $\hat \alpha$ can be made rational by chosing $b$ appropriately. This is however completely obvious,  having in
mind that $\pi^*[T]$ is already rational.  \qed

In case $T$ is the integration over an effective curve, the proof of (2.1) also shows - as already
mentioned

\begin{proposition} Let $\pi: \hat X \to X$ be as in (2.1), $C \subset X$ an effective curve. 
If $[C] \in \Int \sK^*(X),$ then 
there exists an effective curve $\hat C  \subset \hat X$ with $\pi_*(\hat C) = C$ as cycle such that $[\hat C] \in \Int \sK^*(\hat X).$
\end{proposition}

The ``converse" of (2.1) is also true:

\begin{proposition} Let $\pi: \hat X \to X$ be as in (2.1). If $ \sK^*(\hat X) $ contains an interior rational point, so does $\sK^*(X).$
\end{proposition}

\proof Let $\alpha \in \Int \sK^*(\hat X) \cap H^4(\hat X,\bQ)$ and represent it by a positive current $\hat T.$ Let $T = \pi_*(\hat T).$
Then we claim that
$$ [T] \in \Int \sK^*(X) \cap H^4(X,\bQ).$$
It is clear that $[T]$ is rational. Notice that 
$$ \pi^*: H^q(X,\bR) \to H^q(\hat X,\bR) $$
is injective, since $\pi$ is surjective. Hence $\pi^*$ maps ${\overline {\sK(X)}}$ {\it into} ${\overline {\sK(\hat X)}},$ therefore
$$ T(\beta) = \hat T(\pi^*(\beta)) > 0$$
for all $\beta \in {\overline {\sK(X)}}Ê\setminus \{0\},$ proving our claim. \qed

Again we have the following special case:

\begin{proposition} Let $\hat C \subset \hat X$ be an effective curve, $C = \pi_*(\hat C).$ If $[\hat C]Ê\in \Int \sK^*(\hat X),$ then
$[C] \in \Int \sK^*(X)$.  
\end{proposition}

The same reasoning as in (2.3) actually shows (in all dimensions)

\begin{proposition} Let $f: X \to Y$ be a surjective holomorphic map of compact K\"ahler manifolds of positive dimension, let $T$ be a positive closed
current of bidimension $(1,1)$ on $X$ such that $[T] \in \Int \sK^*(X).$ Then $[f_*(T)] \in \Int \sK^*(Y).$ 
\end{proposition}

Putting things together (including (2.1) and applying [OP00], we obtain

\begin{proposition} Let $X$ be a compact K\"ahler threefold, $S$ a smooth compact surface and $f : X \rightharpoonup S$ a dominant meromorphic
map. Suppose that $\Int \sK^*(X) \cap H^4(X,\bQ) \ne \emptyset.$ Then $S$ is projective. 
\end{proposition}

Finally we need also to consider very special singular threefolds. The K\"ahler assumption on $X$ we make is just that $X$ has a K\"ahler
desingularisation. We do not try to define the closure of the K\"ahler cone and its dual but just give the minimum amount of definition we need. 

\begin{definition} {\rm Let $X$ be a normal compact complex space, $\dim X = 3$. 
Let $\pi: \hat X \to X$ be a 
desingularisation with $\hat X$ K\"ahler. Let $C \subset X$ be an effective curve. 
We say that {\it ``$C$ is in the interior of the dual K\"ahler cone
of $X$" iff there exists an effective curve $\hat C \subset \hat X$ such that $\pi_*(\hat C) = C $ and such that $[\hat C] \in \Int \sK^*(\hat X).$ }}
\end{definition}

By (2.2) and (2.4) this definition does not depend on the choice of the K\"ahler desingularisation.  The main point is now

\begin{proposition}Let $X$ be a compact K\"ahler manifold, $Y$ a normal compact complex space and $h: X \to Y$
a Galois covering with Galois group $G$, \'etale in codimension 1. Let $C = \sum a_iC_i $ be an effective curve in $Y$ such that $C_i \not \subset {\rm Sing}Y$ for
all $i.$ Suppose that $C$ is an interior point of the dual K\"ahler cone of $Y.$ Then there 
exists a $G$-stable effective $1$-cycle $C' \subset X$ with $h_*(C') = C$
set-theoretically
such that $[C'] \in \Int \sK^*(X).$ 
\end{proposition} 

\proof First of all notice that $Y$ admits a K\"ahler desingularisation $\pi: \hat Y \to Y. $ 
The existence of $C' $ with $h(C') = C$ is clear; moreover we can make it
$G$-stable by averaging. Hence we may assume $h_*(C') = C$ as cycles.  
Suppose $[C'] \in \partial \sK^*(X).$ Then we find some $\alpha \in {\overline {\sK(X)}} 
\setminus \{0\}$ with $C' \cdot \alpha = 0.$ 
Substituting possibly $\alpha$ by 
$$ \tilde \alpha = \sum_{g \in G} g^*\alpha \in \cK,$$
we may assume from the beginning that $\alpha $ is $G-$invariant, and we still have $C' \cdot \alpha = 0$ since $C'$ is $G$-stable. Thus $\alpha $ is of the form
$$ \alpha = h^*(\beta).$$ We are going to show that this contradicts $[C]$ to be an interior point of $\sK^*(Y).$ By definition and our assumption
that $[C]$ is an interior point, we have an effective curve $\hat C \subset \hat Y$ which is an interior point of $\sK^*(\hat Y)$ and projects to $C.$ On the other
hand we claim
$$ \pi^*(\beta) \in {\overline {\sK(\hat Y)}} \setminus \{0\},  \eqno (1) $$
and
$$ \hat C \cdot \pi^*(\beta) = 0. \eqno (2)$$
This will give a contradiction.\\
To verify (1) we may assume that $\alpha$ is represented by a K\"ahler form $\omega$  
(then take closure). By averaging, we can write $$\omega = h^*(\omega'),$$
where $\omega'$ is a K\"ahler form on $Y_{reg}$ but $\omega'$ extends to all of $Y$ as explained now. Fix a point $y_0 \in Y$ and a
small neighborhood $U$ and consider $V = h^{-1}(U). $ Then on $V$ we can write
$$ \omega = i \partial {\overline {\partial}}Ê\phi$$
with a strictly plurisubharmonic function $\phi.$ By substituting $\phi$ by $\sum_{g \in G} g^*\phi$ and observing $\partial {\overline {\partial}} (g^* \phi) =
g^* \partial {\overline {\partial}} \phi, $ we may assume $\phi $ to be $G$-invariant, hence we can write
$$ \phi = h^*(\phi')$$
with a plurisubharmonic function  $\phi'$ on $U.$ Therefore $$\omega' = i \partial {\overline {\partial }} \phi'$$
on $U_{reg}$ with an extendable function $\phi'$. From this description it immediately follows that $\pi^*(\omega')$ extends
to a semipositive $(1,1)-$form on all of $\hat Y,$ proving (1). \\
(2) We choose a sequence of smooth blow-ups $\rho: Z \to X$ inducing a map $\sigma: Z \to \hat Y$ such that
$$ \pi \circ \sigma = h \circ \rho.$$
Choose an effective curve $\tilde C \subset Z$ with $\rho_*(\tilde C) = C'.$ Then $h_* \rho_*(\tilde C) = C,$ hence $\pi_*\sigma_*(\tilde C) = C$
and thus $\hat C$ and $\sigma_*(\tilde C)$ differ only by cycles supported on the exceptional locus of $\pi.$ Consequently
$$ \hat C \cdot \pi^*(\beta) = \tilde C \cdot \sigma^* \pi^*(\beta) = \tilde C \cdot \rho^*(\alpha) = \hat C \cdot \alpha = 0,$$
proving (2).  \qed

\section{ Interior points and algebraicity} 

Here we investigate the influence of a curve which is an interior point of the dual K\"ahler cone on the algebraicity on the underlying
K\"ahler manifold. 

\begin{theorem} Let $X$ be a compact K\"ahler threefold, $C \subset X$ an effective curve. If $[C] \in \Int \sK^*(X),$ then $a(X) \geq 2$, unless
$X$ is simple non-Kummer. 
\end{theorem}

\proof Suppose that $a(X) \leq 1.$ We have to prove that $X$ is simple non-Kummer. \\
(1) Case $a(X) = 0.$ \\ \sn
If $X$ is not at the same time simple and non-Kummer, then by (1.6), $X$ is uniruled or Kummer. In the first case we have a meromorphic map $f: X \rightharpoonup S$
to a smooth surface $S$ with $a(S) = 0,$ contradicting (2.6). So $X$ is Kummer, hence $X$ is bimeromorphically equivalent to $T/G$ with a torus $T$ and a
finite group $G.$ By (2.2) we may assume that there is a holomorphic bimeromorphic map $f: X \to T/G.$ Then $\dim f(C) = 1.$ This is either
seen by desingularising $T/G$, making $f$ again holomorphic and applying (2.2) or by directly using the arguments of (2.8). Hence
$T$ contains a compact curve which may be assumed $G$-stable. Since $a(T) = 0,$ $T$ has the structure of an elliptic fiber bundle $h: T \to T'$ to a torus $T'$ with 
$a(T') = 0.$ Now $G$ acts also on $T'$ and we have a map $T/G \to T'/G.$ In total we obtain a map $X \to T'/G$ contradicting (2.6). \\

\noindent (2) Case $a(X) = 1.$  \\ \sn
Here we may assume that we have a {\it holomorphic} algebraic reduction $f: X \to B$ to the smooth curve  $B.$ By (2.5), $f(C) = B,$ so $f$
has a multi-section. By (1.6), $X$ is therefore bimeromorphic to $(A \times F)/G =: X',$ where $A$ is a torus or K3 with $a(A) = 0$, $F$ is a smooth
curve and  $G$ a finite group acting diagonally on $A \times F$. Thus we have a holomorphic map 
$$ (A \times F)/G \to A/G$$
and thus a meromorphic dominant map $X \rightharpoonup A/G.$ 
This contradicts Theorem 2.5
It is also possible to avoid the use of 2.5 (and thus of [OP00]) and argue directly by 2.8.
  \qed

We are now addressing the question whether it is possible to have a compact K\"ahler threefold $X$ with $a(X) = 2$ admitting a 
smooth curve $C$ such that $[C] \in 
\Int \sK^*(X)$. In this direction we prove

\begin{proposition} Let $f: X \to S$ be a surjective holomorphic maph connected fibers from the smooth compact K\"ahler 3-fold $X$ to
the smooth projective surface $S.$ Let $C \subset S$ be an irreducible curve with $C^2 > 0.$ Suppose that $X_C = f^{-1}(C)$ is Moishezon and
irreducible. Let $B \subset X_C$ be a general irreducible curve. Then $[B] \in \Int \sK^*(X).$ 
\end {proposition}

A general curve $B \subset X_C$ is understood to be the image of a general hyperplane section $\hat B \subset \hat X_C.$ \\

\proof Suppose $[B] \in \partial \sK^*(X).$ Then there exists a non-zero class $\alpha \in \partial \sK(X),$ such that
$$ \alpha \cdot B = 0.$$
We claim 
$$ \alpha \cdot X_C = 0. \eqno (1) $$
Let $i: X_C \to X$ denote the inclusion and let $\pi: \hat X_C \to X_C$ be a desingularisation; set $g = \pi \circ i.$
Then
$$ g^*(\alpha) \cdot g^*([B]) = 0.$$ 
Now we can write $$g^*([B])= [\hat B] + \sum a_i E_i$$
with $a_i \geq 0$ and $E_i$ the exceptional curves of $\pi.$ 
Hence we conclude that $$ g^*(\alpha) \cdot \hat B = 0.$$ 
Now $g^*(\alpha) \in \overline {\sK(\hat X_C)},$ thus the ampleness of $\hat B$ 
gives $g^*(\alpha) = 0,$ i.e. $\alpha \cdot X_C = 0,$ proving (1). \\
Since $C^2 > 0$ and $C$ is irreducible, some multiple $mC$ is generated by global sections (by a theorem of Zariski, see [Ha70,p.65]).
So $mC$ carries a metric whose curvature form $\omega$ is semipositive and positive outside finitely many curves. By (1) we obtain
$$ \alpha \cdot f^*(\omega) = 0,$$ so that
$$ \alpha \cdot f^*(\omega) \cdot \eta = 0 \eqno(2)$$
for all K\"ahler forms $\eta $ on $X.$ 
In order to exploit this, we represent $\alpha$ by a positive closed $(1,1)$-current $T.$ Let $S_0$ be the maximal open
subset of $S$ over
which $f$ is a submersion and let $X_0 = f^{-1}(S_0).$ Then a small local cacluation in $X_0$ shows that $f^*(\omega) \wedge \eta$
is a strictly positive $(2,2)$-form on $X_0.$ By (2) we have
$$ T(f^*(\omega) \wedge \eta ) = 0,$$
hence $\dim {\rm supp}(T) \subset X \setminus X_0.$  Since $\dim (X \setminus X_0) = 2,$
Siu's well-known theorem says that 
$$ T = \sum \lambda_i Z_i,$$
where the $Z_i$ are irreducible components of $X \setminus X_0$, where $\lambda_i > 0$ and where the right hand side of course
means integration over the cycle. We fix a positive integer $m$ such that $mC$ is generated by sections and therefore defines a 
map $g: S \to S'$ to a normal projective surface $S'.$ Now let $C' \subset \vert mC \vert $ be a general member. 
Then (2) yields
$$ \sum \lambda_i Z_i \cdot X_{C'} = 0.$$
Since the support of both divisors $\sum \lambda_i Z_i$ and $X_{C'}$ do not have common components, it follows that $Z_i \cap X_{C'} =
\emptyset $ for all $i.$ Now either $\dim f(Z_i) = 0$ or $ f(Z_i) \cdot C' = 0, $ whence $f(Z_i)$ has to be contracted by $g.$
So
$$\dim (g \circ f)(Z_i) = 0 $$
for all $i.$ This is incompatible with $\alpha = \sum \lambda_i [Z_i] \in \overline {\sK (X)}, $ and thus producing a contradiction. 
In fact, we decompose $\sum  \lambda_i Z_i$ into connected components and write with the appropriate multiplicities:
$$ \sum \lambda_i Z_i = \sum \mu_j W_j.$$
So the $W_j$ are pairwise disjoint and are supported in fibers of $g \circ f$ so that clearly all $W_j$ are not nef (recall that the
general fiber of $g \circ f$ has dimension $1$). So $\alpha \not \in \overline {\sK (X)}.$  \qed

With the same proof we can also deal with the case that $X_C$ is irreducible:

\begin{proposition} If in 3.2 $X_C$ is reducible, say $X_C = A_1 \cup \ldots \cup A_r,$ then we take $B_i \subset A_i$ general and
put $B = \sum B_i.$ Again $[B] \in \Int \sK^*(X).$ 
\end{proposition}

\begin{corollary} Let $X$ be a smooth compact K\"ahler 3-fold with $a(X) = 2$ with holomorphic algebraic reduction $f: X \to S$ to a
smooth projective surface $S.$ Let
$$ \Delta = \overline { \{s  \in S \vert X_s {\rm \ singular, \ not  \ multiple \ elliptic \ } \}}.$$
Suppose $\Delta$ contains an irreducible component $C$ with $C^2 > 0.$ Then $X$ carries a (possibly reducible) curve $B$ such that 
$[B] \in \Int \sK^*(X).$ If $f^{-1} (C)$ is irreducible, i.e. the general fiber over $C$ is irreducible (rational), then $B$ can be
taken irreducible.
\end{corollary}

\proof We only have to show that $X_C = f^{-1}(C)$ is Moishezon. This however is clear since $X_C$ is covered by rational curves,
the fibers of $f$ over $C.$ \qed  

\begin{remark} {\rm (1) It remains to explicitly construct a K\"ahler 3-fold $X$ with $a(X) = 2$ (with holomorphic algebraic reduction; this
can always be achieved by blowing up $X$) fulfilling the requirements of (3.4). It seems unimaginable that such a 3-fold should not exist;
non-existence would simply mean that any component $C$ of $\Delta$ (as in (3.4)) has $C^2 \leq 0.$ On the other hand, an explicit construction
is not obvious. Of course there are projective elliptic 3-folds with a curve $C \subset  \Delta$ and $C^2 > 0.$ So one might try to 
deform $X$ in a non-algebraic K\"ahler 3-fold keeping the curve $C.$ \\
(2) It is easy to construct elliptic K\"ahler 3-folds $f: X \to S$ of algebraic dimension 2 admitting a curve $C \subset  \Delta$ such that
$C^2 = 0.$ In fact, let $f_1 : S_1 \to B$ be an elliptic K\"ahler surface with $a(S) = 1 $ and $f_2 : S_2 \to B$ a projective surface. Suppose
that $f_1$ has some non-elliptic fibers
and that the singular loci of $f_1 $ and $f_2$ in $B$ are disjoint. 
Then
$$ X = S_1 \times_B S_2$$
is a smooth K\"ahler 3-fold with $a(X) = 2$ and algebraic reduction $f: X \to S_2.$ Now $\Delta_f = \Delta_{f_1} \times_B S_2 $ so that
$\Delta_f$ contains of fibers of $f_2.$ \\
(3) Here is a possible way how to construct $g: X \to S$ as required in (3.4). The procedure is as follows: we start with a fibration
$f: X \to C$ from the compact K\"ahler 3-fold $X$ with $a(X) = 2$ to the smooth curve $C.$ Suppose that $a(X_c) = 1$ for the general fiber $X_c$
of $f$ such that the algebraic reduction has non-elliptic singular fibers. Then we can form the relative algebraic reduction $X \rightharpoonup S \to C.$ 
Suppose that $S$ is smooth and $f: X \rightharpoonup S$
is actually holomorphic. Then $f$ is an elliptic fibration and $S$ is projective. Moreover ``in general'' $\Delta$ should dominate $C$
This gives hope to find $C \subset \Delta $ with $C^2 > 0.$ The difficulty is to find a starting fibration $g: X \to C;$ the possible most 
natural choice would be require $X_c$ to be a Kummer surface and that $g$ is a submersion. Then the above construction will work and therefore
we are reduced to the following}
\end{remark}

\begin{problem} Does there exist a compact K\"ahler threefold $Z$ with a submersion $h$ to a curve $C$ such that the general fiber is a torus
of algebraic dimension 1? 
\end{problem}

It is not so difficult to construct a non-K\"ahler submersion $Z \to C$ such that the general fiber is a torus of algebraic dimension 1,
but the K\"ahler property seems difficult to achieve.

\bigskip 

Here is a mild restriction for K\"ahler 3-folds $X$ containing a curve $C$ which is an interior point of the dual K\"ahler
cone.

\begin{proposition} Let $X$ be a smooth compact K\"ahler 3-fold with algebraic reduction $f: X \rightharpoonup S$ to the smooth 
projective surface $S.$ Suppose $\Int \sK(X)$ contains a (possibly reducible) curve. Then $q(X) = q(S).$ 
\end{proposition}

\proof Assume $q(X) > q(S).$
 By considering a holomorphic model for $f$ and appyling the Leary spectral sequence, we see immdiately that $q(X) = q(S) + 1.$ Let $q = q(S).$ 
If $\alpha : X \to \alpha (X) \subset A_{q+1}$ and $\beta: S \to \beta (S) \subset B_q$ denote the Albanese maps, then we 
conclude that $\alpha (X)$ cannot be projective, otherwise we would have $q(S) = q(X).$ Let $\gamma: \alpha (X) \to \beta (Y)$
be the induced map. Since $\beta (Y) $ is projective, $\gamma$ is not generically finite. It is important to notice that 
$\dim \gamma \alpha (C) = 1$ by (2.4).

 We are thus left with the following cases: \\
(a)$ \dim \alpha (X) = 2.$ Then $\alpha (C)$ is a multi-section of $\gamma,$ making $\alpha (X)$  projective. \\
(b) $\dim \alpha (X) = 3.$ Then necessarily $\dim \beta (Y) = 2$ and now $\gamma$ has a partial multi-section $\alpha (X).$
But $\gamma $ is the reduction of the canonical projection $A_{q+1} \to B_q $ which (up to finite \'etale cover) is an
elliptic fiber bundle without transverse curves; contradiction. \qed

\section{ Interior points and the normal bundle}

In this section we investigate curves in K\"ahler manifolds, their normal bundles and the connection with the dual K\"ahler cone. Here it is not always necessary
to restrict to low dimensions.

\begin{conjecture} Let $X$ a compact K\"ahler manifold and $C \subset X$ a smooth curve with ample normal bundle.
Then $[C] \in \Int \sK^*(X).$ \end{conjecture}

As a direct consequence of Theorem 3.1 in dimension three, we are lead to 

\begin{conjecture} Let $X$ be a compact K\"ahler manifold, $Y \subset X$ a positive-dimensional compact submanifold  with ample normal bundle. Then $a(X) \geq \dim Y +1. $
\end{conjecture}

In our feeling, (4.2) seems more accessible than (4.1). Notice that (4.2) holds for hypersurfaces. 
\smallskip \noindent

In this section we verify (4.1) and (4.2) in special cases. First we notice that (4.1) and (4.2) hold in the
surface case, even without K\"ahler assumption.

\begin{proposition} Let $X$ be a smooth compact surface, $C \subset X$ an effective curve with ample normal bundle.
Then $X$ is projective and $[C] \in \Int \ \sK^*(X).$ \end{proposition}

\proof The first assertion is classical using Riemann-Roch and the fact that Moishezon surfaces are projective.\\
The second assertion is clear if $C$ is an ample divisor on $X.$
If $C$ is not ample, $C$ is big and nef, in particular $C^2 > 0.$ Now assume that $[C]$ 
is not an interior point of $\sK^*(X).$
Then we find $0 \ne \alpha \in \overline {\sK(X)}$ such that
$$ C \cdot \alpha = 0.$$
Hence the Hodge Index Theorem for $H^{1,1}$ yields $\alpha \equiv mC$ with a positive real number
$m.$ So $C \cdot \alpha = m C^2 = 0,$ contradiction. 

\qed

Now we prove a special case of conjecture 4.1. 

\begin{theorem} (1) Let $C \subset X$ be an effective curve in the smooth projective threefold $X.$ Suppose that there is
a smooth ample surface $S \subset X$ with $C \subset S$ such that $N_{C\vert S}$ is ample (e.g. $C$ is a complete intersection
$C = S \cap S'$ with ample $S'$). Then $[C] \in \Int \  \sK^*(X).$ \\
(2) Let $C$ be a smooth complete intersection curve of hyperplane sections in the projective manifold $X.$ Then $[C] \in \Int \  \sK^*(X).$
\end{theorem}

\proof (1) The inclusion $i: S \to X$ yields a map $i^* : {\overline {\sK (X)}} \to 
{\overline {\sK (S)}}$. Since $S$ is ample, $i^*$ is
injective.  
Take $\alpha \in \cK$ and represent it by a form $\omega.$ Then $\int_C \omega = \int_C i^*(\omega) ,$ and since
$[i^*(\omega)] \in \overline {\sK(S)},$ we conclude from (4.3) that $ \alpha \cdot C = \int_C \omega > 0.$ \\
(2) follows now by induction. \qed

Concerning Conjecture 4.2 we prove

\begin{theorem} Let $C \subset X$ be a smooth curve with ample normal bundle $N_C$ in the 
compact K\"ahler manifold $X.$ Assume moreover that
$N_C \otimes T_C$ is ample.
Then $X$ is projective. \end{theorem}

\proof We check the projectivity of $X$ by showing
$$ H^2(X,\sO_X) = H^0(X,\Omega^2_X) = 0. \eqno (*)$$
 
This is shown by adopting the method of [PSS99, 2.1]: by power series expansion, one has the 
inequality
$$ h^0(X,\Omega^2_X) \leq \sum_{k=0}^{\infty} h^0(C,S^kN_C^* \otimes \Omega^2_X \vert C).  \eqno (**)$$
Taking $\bigwedge^2 $ of the sequence
$$ 0 \to N^*_C \to \Omega^1_X \vert C \to \Omega^1_C \to 0, $$
we obtain a sequence
$$ 0 \to \bigwedge^2N^*_C \to \Omega^2_X \vert C \to \Omega^1_C \otimes N^*_C \to 0.$$ 
Tensoring with $S^kN^*_C$ for $k \geq 0$, taking cohomology and having in mind our ampleness assumptions, we obtain the vanishing
$$ H^0(C,S^kN^*_C \otimes \Omega^1 \otimes N^*_C) = 0$$
for $k \geq 0$
from which (*) follows by virtue of (**). \qed

We remark that in most cases in Theorem 4.5 it is easily seen that $[C] \in {\rm Int}\sK^*(X).$

\begin{corollary} Let $X$ be a compact K\"ahler manifold, $C \subset X$ a smooth curve with ample normal bundle. If $g(C) \leq 1,$ then
$X$ is projective.
\end{corollary}

Of course one can formulate more general versions of (4.5) e.g. for reduced locally complete intersection curves; we leave that to the reader.
We will prove more on the structure of $X$ (in case $C$ is elliptic and $\dim X = 3$) in sect.5.
Observe that (4.5) and (4.6) hold also in the following more general context: it is sufficient to assume that $X$ is in class ${\cal C}$, i.e. bimeromorphic
to a K\"ahler manifold. Then the conclusion is that $X$ is Moishezon. 

\begin{theorem} Let $X$ be a compact K\"ahler threefold and $C \subset X$ a smooth curve. Suppose that the normal bundle $N_C$ is ample. \\
(1) If $\kappa (X) \leq 0,$ then $X$ is projective or $X$ is simple non-Kummer. \\
(2) $a(X) \ne 1.$ 
\end{theorem}

\proof (I) We first show the following general statement. Suppose $f: X \to Y$ is a surjective holomorphic map,
then $\dim f(C) = 1.$ \\
In fact, suppose $C \subset f^{-1}(y)$ for some $y \in Y.$ Let $F$ be the complex-analytic fiber over $y.$ Now choose
$k$ maximal such that the $k-$th infinitesimal neighborhood $C_k$ is still a subspace of $F.$ Then
we obtain an exact sequence of conormal sheaves
$$ N^*_{F/X}\vert C_k \to N^*_{C_k/X} \to N^*_{C_k/F} \to 0.$$ 
Restricting to $C$ and observing that
$$ N^*_{F/X}\vert C \to N^*_{C_k/X} \vert C $$
does not vanish by the choice of $k,$ the spannedness of $N^*_{F/X}$ provides a non-zero section of $ N^*_{C_k/F} \vert C.$
Thus $S^kN^*_{C/X}$ has a non-zero section, contradicting the ampleness of $N_C.$

\bn (II) We next show that $a(X) \ne 1.$ We first claim that every resolution of indeterminacies of
the algebraic reduction $f: X \rightharpoonup B$ admits a multi-section.
In fact, otherwise $f$ itself is holomorphic and would contract $C$ contradicting (I). Then by (1.6), 
$X$ is bimeromorphically of the form $(A \times F)/G,$ with $a(A) = 0$ and $F$ a curve. So $A$ is K3 or a torus. 
In case $A$ is K3, we subsitute $A$ by the blow-down of all $(-2)-$curves. So we may assume that $A$ does not
contain any curve at the expense that $A$ may be singular. Therefore the induced map $$h: X \rightharpoonup A/G$$
is actually holomorphic. In particular $h$ contracts $C$ which again contradicts (I).

\bn (III) Now suppose $\kappa (X) = - \infty$ and that $X$ is not both simple non-Kummer. \\
 Suppose moreover $a(X) = 0.$ 
In that case we have by (1.6) a holomorphic map $f: X \to S$ to a normal surface $S$ with $a(S) = 0,$
the general fiber of $f$ being a smooth rational curve. Arguing as in (II), we may assume that $S$ does not
contain any curve. So $f(C)$ is a point, contradicting (I). So $a(X) \ne 0.$ \\
 By (II) it therefore remains to exclude $a(X) = 2.$
So suppose $a(X) = 2$ and let $f: X \rightharpoonup S$ be an algebraic reduction. By [CP00], $X$ is uniruled, so
we can form the rational quotient $h: X \rightharpoonup Z$ of a covering family of rational curves. 
The fibers of $h$  being rationally connected,
we conclude that ${\rm dim}Z > 1,$ otherwise $X$ would be projective. Therefore $\dim Z = 2$ and since $a(X) = 2,$
we must have $a(Z) = 1.$ 
Using the algebraic reduction $Z \to B$ we obtain a meromorphic map $g: X \rightharpoonup B.$ The general fiber of $g$ is
bimeromorphically a $\bP_1$-bundle over an elliptic curve, hence algebraic. Again a holomorphic model of $g$ cannot have
a multi-section by (1.2), hence $g$
must be holomorphic and again provides a contradiction. \qed

\bn (IV) We are left with $\kappa (X) = 0$ and suppose $X$ is not both simple non-Kummer.
By [CP00], Theorem 8.1, $X$ has a birational model $X'$ which admits a finite cover, \'etale in codimension 1, say $\tilde X$, such that $\tilde X$ is either
a torus or a product of an elliptic curve with a K3 surface. \\
(a) $a(X) = 0.$ Then $X$ is bimeromorphically $T/G$ with $T$ a torus. If $T$ has no curves, then we actually have a
holomorphic map $h: X \to T/G$ and conclude by (I). If $T$ admits a curve, we have an elliptic fiber bundle structure
$g: T \to T'$ to a 2-dimensional torus $T'$ with $a(T') = 0$ and every curve in $T$ is a fiber of $g.$ Then there is 
an induced map $T/G \to T'/G,$ and $T'/G$ has no curves, so that the meromorphic map $X \rightharpoonup T'/G$ is
actually holomorphic. We conclude by (I). \\
(b) Since $a(X) \ne 1 $ by (II), we are left with $a(X) = 2.$ \\
(b.1) Suppose first $\tilde X = W \times E$ with $W$ a K3-surface, $a(W) = 1$ and $E$ elliptic. With $X' = \tilde X/G$, we note
that the $G-$action  on $\tilde X$ is diagonal, hence we have by composing with the algebraic reduction $h: W/G \to C$ 
a meromorphic map $$g: X \rightharpoonup C$$
to a rational curve $C.$ Since $h$ has no multi-sections, any multi-section of a resolution of $X \to W/G$ must be
mapped to a fiber of $h.$ Hence $g$ is almost holomorphic, hence holomorphic. Since $\dim g(C) = 0,$ we conclude by (I).\\
(b.2) Suppose finally that $\tilde X$ is a torus. Then the algebraic reduction of $\tilde X$ provides an elliptic bundle
structure $\tilde X \to \tilde B$ to an abelian surface without multi-sections. We obtain a map $\tilde X/G \to \tilde B/G$
without multi-sections and $X \sim \tilde X/G.$ So the induced map $X \rightharpoonup \tilde B/G$ must be
holomorphic and will contract $C$, contradicting (I). \qed

\begin{corollary} Let $X$ be a compact K\"ahler threefold, $C \subset X$ a smooth curve. Assume that the normal
bundle $N_C$ is ample. Then $a(X) \geq 2$ unless $X$ is simple non-Kummer.
\end{corollary}

\begin{corollary} Let $C$ be a smooth curve in the compact K\"ahler threefold with ample normal bundle. If
$K_X \cdot C < 0,$ then $X$ is projective unless $X$ is simple non-Kummer with $\kappa (X) = - \infty.$
\end{corollary}

\proof It is sufficient to prove $\kappa (X) = - \infty;$
then we apply (4.7).
This follows essentially from [PSS99,2.1] and is parallel to the argument in (4.5).
Namely, we claim that for all $t \in \bN $ and all $k \in \bN$ the following vanishing holds
$$ H^0(C,S^kN^*_C \otimes K_X^t \vert C) = 0. \eqno (*)$$

This is clear since $N_C$ is ample and $K_X \vert C$ is negative by adjunction. 
On the other hand power series expansion gives

$$ h^0(X,K_X^t) \leq \sum _{i=0}^{\infty} h^0(C,S^kN^*_C \otimes K_X^t \vert C). $$
Thus (*) gives $\kappa (X) =  - \infty.$ \qed

For the cases $\kappa (X) = 1,2$ we need a stronger assumption than just the existence of the family $(C_t).$

\begin{theorem} Let $X$ be a compact K\"ahler threefold, $C \subset X$ a smooth curve with ample normal bundle.
Suppose that $C$ moves in a family $(C_t)$ that {\bf covers} $X.$ Then $X$ is projective.
\end{theorem}

\proof By (4.7) we may assume $a(X) \geq 2.$
So we have a meromorphic elliptic fibration $f: X \rightharpoonup S$ to a projective
surface $S.$ Since $g(C_t) \geq 2,$ we conclude $\dim f(C_t) = 1$ for general $t$ and therefore the $C_t$ make $X$ algebraically connected, thus
$X$ is projective, contradiction. \qed

\begin{example} {\rm In general curves in the interior of the dual K\"ahler cone of course are far from having ample normal bundle.
For a trivial example take a threefold $X$ with $b_2(X) = 1$ containing a $(-1,-1)-$curve, i.e. a smooth rational curve C with
normal bundle $N_C = \sO(-1) \oplus \sO(-1).$ } 
\end{example}

Another interesting class of curves are the connecting curves:

\begin{definition} A family $(C_t)$ of $1-$cyycles (say with general $C_t$ reduced) is a connecting family if and only if any general points $x,y \in X$ 
can be joined by a chain of curves of type $C_t.$
\end{definition}

By (1.2), any variety $X$ carrying a connecting family of $1-$cycles is projective.

\begin{theorem}  Let $X$ be a compact K\"ahler manifold, $(C_t)$ a connecting family, then $[C_t] \in \Int \sK^*(X)$. 
\end{theorem}

\proof [WK01].  \qed

\begin{remark} {\rm Having in mind the results of sect.3, it is likely that there exists a compact K\"ahler 3-fold $X$ with
$a(X) = 2$ carrying a smooth or at least a locally complete intersection curve with ample normal bundle. The difficulty of course
is that in 3.4 the surface $X_C$ need not be projective. \\
Notice also the following.
If $C \subset X$ is a smooth curve with ample normal bundle in a projective manifold, then ${\cal M}(\hat C)$ is a finite-dimensional
${\cal M}(X)-$ vector space, where $\hat C$ is the formal completion of $X$ along $C$ and ${\cal M}$ is the sheaf of meromorphic functions. 
One says that $C$ is $G2$ in $X,$ see [Ha70]. A counterexample as above with ample normal bundle would indicate that this property does no 
longer hold in the K\"ahler setting. In fact, it seems likely that the ampleness of the normal bundle forces the field  of formal meromorphic
functions along $C$ should have transcendence degree $3$ over the complex numbers. }

\end{remark}

For some new results concerning Conjecture 4.1 we also refer to [BM04].

\section{Structure of projective threefolds containing a smooth elliptic curve with ample normal bundle}

Let $X$ be a smooth projective threefold containing a smooth curve $C$ with ample normal bundle $N_C.$ If $C$ is rational, 
then $X$ is rationally connected [KoMiMo92].
In this section we consider the case that $C$ is elliptic and we shall fix this situation unless otherwise stated. By [PSS99,2.1], see (4.5), we have $\kappa (X) = - \infty,$ hence $X$ is uniruled. Applying the minimal model
program, we find a birational rational map $f: X \rightharpoonup X'$, given by a sequence of birational contractions and flips, and
a contraction $g: X' \to Y$ with $\dim Y \leq 2.$ Notice that in case $\dim Y = 2,$ $Y$ has only rational singularities and also that
$\pi_1(X) = \pi_1(X') = \pi_1(Y)$ (see e.g. [Ko95]). We have the following structure theorem.

\begin{theorem} (1) The irregularity $q(X) \leq 1.$ \\
(2) If $q(X) = 0$  then $X$ is rationally connected.
In particular we always have $\pi_1(X) = 0$. \\
(3) If $q(X) = 1$, then $C$ is an \'etale multi-section of the Albanese map $\alpha: X \to A$ and $\alpha $ factor over the Albanese $\beta: Y \to A.$ The general
fiber of $\alpha $ is a rational surface.
Moreover either $\beta = {\rm id},$ or $\dim Y = 2,$ $\kappa (Y) = - \infty$ and $Y$ is a generic $\bP_1-$bundle over $A.$ In particular
$\pi_1(X) = \zed^2;$  more precisely 
$$ \alpha_*: \pi_1(X) \to \pi_1(A) = \zed^2 $$
is an isomorphism and ${\rm Im}(\pi_1(C) \to \pi_1(X))$ has finite image. 
\end{theorem}

\proof (1) It is known in general that ${\rm Alb}(Z) \to {\rm Alb}(X)$ is surjective, if $Z \subset X$
is a submanifold with ample normal bundle in $X$ (see e.g. [Ha70,p.116]). Applying this to $Z = C,$
we get $q(X) \leq q(C) = 1.$ \\
(2) By [KoMiMo92], it is sufficient to show that
$$ H^0(X,S^t\Omega^2_X) = 0  \eqno (*)$$
for all $t \geq 1$ (actually $t = 2$ suffices). We verify this by the same method as in (4.5). In fact, 
(*) will follow from 
$$ H^0(C,S^kN^*_C \otimes S^t\Omega^2_X \vert C) = 0  \eqno (**)$$
for all $k \geq 0 $ and all $t \geq 1.$ 
In order to verify (*), we use the sequence 
$$ 0 \to N_C^* \to \Omega^1_X \vert C \to \Omega^1_C = \sO_C \to 0 $$
and take $\bigwedge^2$ to obtain 
$$ 0 \to  {\rm det}N^*_C \to \Omega^2_X \vert C \to N^*_C \otimes \Omega^1_C \to 0.$$
Hence $\Omega^2_X \vert C$ is a negative vector bundle and thus (**) is clear. \\
(3) Suppose now $q(X) = 1.$ Then $C$, having ample normal bundle, is not contracted by $\alpha$ and therefore it is an \'etale multi-section
of $\alpha.$ The existence of $\beta $ and the factorisation property are clear. If $\dim Y = 1,$ 
then we must have $Y = A$ and $\pi_1(X) =
\pi_1(A) = \zed^2.$ So let $\dim Y = 2.$ Notice that (**) from (2) was independent of $q(X) = 0,$ 
so it holds also in our context.
This implies $\kappa (Y) = - \infty,$ whence our claim. \qed

\begin{re} {\rm Notice that the proof of (5.1) yields actually the following. Let $X$ be a projective manifold, $C \subset X$ an elliptic curve with
ample normal bundle. Then 
$$ H^0(X,S^t\Omega^k_X) = 0$$
for all $t \geq 1$ and all $k \geq 2.$ In particular, $X$ does not admit any rational map to a variety $Y$ with $\dim Y \geq 2$ and $\kappa (Y) \geq 0.$}
\end{re} 

In this context we notice the following very interesting result of Napier-Ramachandran [NR98]:

\begin{theorem}{\it (Napier-Ramachandran)} Let $X$ be a projective manifold, $Y \subset X$ a submanifold of positive dimension with ample normal bundle $N_Y.$Then
$$ {\rm Im}(\pi_1(Y) \to \pi_1(X)) $$
has finite index. 
\end{theorem}

Observe that the surjectivity ${\rm Alb}(Y) \to {\rm Alb}(X)$, used already in the proof of (5.1), means that 
$$ {\rm Im}(H_1(Y,\zed ) \to H_1(X,\zed )) $$
has finite index. This is the abelianized version of (5.3). If the Shafarevich conjecture holds 
(``the universal cover of a projective (compact K\"ahler)
manifold is holomorphically convex"), then by Koll\'ar, see e.g. [Ko95], the so-called Shafarevich map 
$sh: X \to Sh(X)$ exists which 
contracts exactly those subvarieties $Z \subset X$ for which ${\rm Im}(\pi_1(Z) \to \pi_1(X))$ has finite index. Suppose $\pi_1(X) $ is not finite, i.e.$Sh(X)$ 
is not a point - in that case (5.3) is anyway trivial. Since $N_Y$ is ample, it follows
that $\dim sh(Y) \ne 0,$ hence ${\rm Im}(\pi_1(Y) \to \pi_1(X)) $ is at least infinite. Of course we used the ampleness here only in a very
weak form, namely to conclude that $\dim sh(Y) > 0.$ Actually $\dim sh(Y) = \dim Y.$

\vskip 20pt
Keiji Oguiso \\
Department of mathematics, University of Tokyo,\\
Komaba,Meguro,153-8914, Japan \\
e-mail: oguiso@ms.u-tokyo.ac.jp
\par
\vskip 10pt
Thomas Peternell \\
Mathematisches Institut;  Universit\"at Bayreuth \\
D-95440 Bayreuth, Germany \\
e-mail: thomas.peternell@uni-bayreuth.de

\end{document}